\documentclass[11pt]{amsart}

\usepackage{graphicx}
\usepackage{wrapfig}
\usepackage{amssymb}
\usepackage{times}
\usepackage[top=20mm, bottom=20mm, left=20mm, right=20mm]{geometry}
\usepackage{color}
\usepackage{bm}
\usepackage{multirow, bigdelim}
\usepackage{algorithmic}

\newcommand{\mtx}[1]{\bm{\mathsf{#1}}}
\renewcommand{\vec}[1]{\bm{\mathsf{#1}}}

\theoremstyle{definition}
\newtheorem{remark}{Remark}

\newcommand{\lsp}{\vspace{3mm}}

\begin{document}

\begin{center}
\textbf{\large{Efficient nuclear norm approximation via the randomized UTV algorithm}}

\vspace{5mm}

N.~Heavner\footnotemark[1],
P.G.~Martinsson\footnotemark[1]

\vspace{5mm}

\begin{minipage}{140mm}
\textbf{Abstract:} The recently introduced algorithm \texttt{randUTV} provides a highly efficient technique for computing accurate approximations to all the singular values of a given matrix $\mtx{A}$. The original version of \texttt{randUTV} was designed to compute a full factorization of the matrix in the form $\mtx{A} = \mtx{U}\mtx{T}\mtx{V}^{*}$ where $\mtx{U}$ and $\mtx{V}$ are orthogonal matrices, and $\mtx{T}$ is upper triangular. The estimates to the singular values of $\mtx{A}$ appear along the diagonal of $\mtx{T}$. This manuscript describes how the \texttt{randUTV} algorithm can be modified when the only quantity of interest being sought is the vector of approximate singular values. The resulting method is particularly effective for computing the nuclear norm of $\mtx{A}$, or more generally, other Schatten-$p$ norms. The report also describes how to compute an estimate of the errors incurred, at essentially negligible cost.
\end{minipage}

\end{center}

\footnotetext[1]{Department of Applied Mathematics, University of Colorado at Boulder, 526 UCB, Boulder, CO 80309-0526, USA}

\lsp

% =============================================================================
\section{Overview}
% =============================================================================

This note describes an efficient algorithm for computing an accurate estimate for the nuclear norm $\| \cdot \|_*$ of a given matrix $\mtx{A} \in \mathbb{R}^{m \times n}$. The nuclear norm has recently found uses in numerical optimization, beginning with the introduction by Fazel et al.~ in \cite{fazel2002matrix,fazel2001rank} of its use as an effective heuristic for solving the rank minimization problem
\begin{align*}
\text{minimize} &\phantom{=} \text{rank } \mtx{X} \\
\text{subject to} &\phantom{=} \mtx{X} \in \mathcal{C}
\end{align*}
where $\mtx{X} \in \mathbb{R}^{m \times n}$ is the decision variable and $\mathcal{C}$ is some given convex constraint set. Recht et al.~later proved in \cite{recht2010guaranteed} that in certain cases, minimizing the nuclear norm also yields the theoretical solution to the corresponding rank minimization problem, solidifying the validity of the heuristic. The algorithm discussed in this note, \texttt{randNN}, may be used in a line search to choose the step size for nuclear norm minimization algorithms such as projected subgradient methods \cite{recht2010guaranteed} or mirror descent \cite{nemirovskii1983problem,beck2003mirror}.

The recently proposed algorithm \texttt{randUTV} \cite{martinsson2017randutv} is designed to compute, given a  matrix $\mtx{A} \in \mathbb{R}^{m \times n}$, a factorization of the form
\begin{equation}
\label{eq:UTVintro}
\begin{array}{ccccccccccc}
\mtx{A} &=& \mtx{U} & \mtx{T} & \mtx{V}^{*}.\\
m\times n && m\times m & m\times n & n\times n
\end{array}
\end{equation}
In (\ref{eq:UTVintro}), the matrices $\mtx{U}$ and $\mtx{V}$ are orthogonal, $\mtx{T}$ is upper triangular, and the diagonal entries of $\mtx{T}$ are good approximations to the singular values of $\mtx{A}$. This note describes a modified version of \texttt{randUTV} which is of particular use in the sub-problem of conducting a line search during nuclear norm minimization. Specifically, we re-derive the \texttt{randUTV} algorithm from the perspective of nuclear norm estimation. In this case, since only the singular value estimates for $\mtx{A}$, rather than the entire matrix factorization, are desired, several steps of \texttt{randUTV} may be omitted or modified to yield \texttt{randNN}, which sees modest acceleration over \texttt{randUTV} and major acceleration over a full (SVD) computation of the singular values. We also mention a rough upper bound on the accuracy of the computed singular values that may computed as a part of \texttt{randNN} at very little extra cost from the middle matrix $\mtx{T}$.

The structure of this note is as follows. In Section \ref{sec:prelims}, we review the notation used throughout the remainder. We derive the algorithm \texttt{randNN} in \ref{sec:alg}, and in \ref{sec:err_bds} we discuss an error bound for the resulting approximations to the singular values. Finally, Section \ref{sec:numerics} contains numerical experiments exploring the accuracy of the singular values estimated by \texttt{randNN} and the performance of the algorithm compared to the industry standard.

% =============================================================================
\section{Preliminaries}
\label{sec:prelims}
% =============================================================================

In this note, we use the notation $\mtx{A} \in \mathbb{R}^{m \times n}$ to denote a matrix $\mtx{A}$ of dimension $m \times n$ with real entries. $A_{i,j}$ denotes the $i$,$j$th entry of $\mtx{A}$. We use the notation of Golub and Van Loan \cite{golub2012matrix} to specify submatrices: If $\mtx{A}$ is an $m\times n$ matrix, and $I = [i_{1},\,i_{2},\,\dots,\,i_{k}]$ and $J = [j_{1},\,j_{2},\,\dots,\,j_{\ell}]$ are index vectors, then $\mtx{A}(I,J)$ denotes the corresponding $k \times \ell$ submatrix. We let $\mtx{A}(I,:)$ denote the matrix $\mtx{A}(I,[1,\,2,\,\dots,\,n])$, and define $\mtx{A}(:,J)$ analogously. The transpose of a matrix $\mtx{A}$ is denoted $\mtx{A}^*$, and we use $\sigma_i(\mtx{A})$ to reference the $i$th leading singular value of $\mtx{A}$.

We measure matrices with either the spectral norm $\| \cdot \|,$ nuclear norm $\| \cdot \|_*$, or Frobenius norm $\| \cdot \|_F$, with definitions given by
\[
\| \mtx{A} \| = \sup_{\| \vec{x} \| = 1} \| \mtx{A} \vec{x} \|_2 = \sigma_1(\mtx{A}), \quad
\| \mtx{A} \|_* =  \sum_{i=1}^{\min(m,n)} \sigma_i(\mtx{A}), \quad
\| \mtx{A} \|_F = \left (  \sum_{i=1}^m \sum_{j=1}^{n} A_{i,j}^2  \right )^{1/2}
\]
(where $\| \cdot \|_2$ is the Euclidean 2-norm for vectors in $\mathbb{R}^n$).

% =============================================================================
\section{Description of the algorithm}
\label{sec:alg}
% =============================================================================

% -----------------------------------------------------------------------------
\subsection{Approach}
% -----------------------------------------------------------------------------

Consider a matrix $\mtx{A} \in \mathbb{R}^{m \times n}$. Now partition $\mtx{A}$ as

$$\begin{tabular}{ccccccc}
 & & $b$ & \ldelim[{2}{3mm} &  $\mtx{A}_{11}$ & $\mtx{A}_{12}$ & \rdelim]{2}{3mm}\\
$\mtx{A}$ & = & $m - b$ & & $\mtx{A}_{21}$ & $\mtx{A}_{22}$ &  \\
 & & & & $b$ & $n - b$ &
\end{tabular}.$$

The construction of the algorithm will be guided by several fundamental observations:

\begin{itemize}
\item If $\mtx{A}_{12} = \mtx{0}$ and $\mtx{A}_{21} = \mtx{0}$, then $\| \mtx{A} \|_* = \| \mtx{A}_{11} \|_* + \| \mtx{A}_{22} \|_*$.
\item Multiplying $\mtx{A}$ by an orthogonal matrix does not change the nuclear norm.
\item Consider the matrix $\tilde{\mtx{A}} = \mtx{U}^* \mtx{A} \mtx{V}$, where $\mtx{U}$ and $\mtx{V}$ are orthogonal matrices whose first $b$ columns span the same space as the leading $b$ left and right singular vectors of $\mtx{A}$, respectively. Then giving $\tilde{\mtx{A}}$ the same partition as $\mtx{A}$, we have $\tilde{\mtx{A}}_{12} = \mtx{0}$ and $\tilde{\mtx{A}}_{21} = \mtx{0}$. For this reason, we will call such a $\mtx{U}$ and $\mtx{V}$  the ``optimal'' rotation matrices for the problem.
\end{itemize}

These remarks lay out a path for estimating $\| \mtx{A} \|_*$. We will choose $b$ to be a relatively small block size, say 50 or 100, and will compute an orthonormal matrices $\mtx{U}$ and $\mtx{V}$ whose first $b$ columns \textit{approximately} span the same space as the respective leading $b$ left and right singular vectors of $\mtx{A}$. We will then form $\mtx{A}^{(1)} = \mtx{U}^* \mtx{A} \mtx{V}$, after which we will have $\mtx{A}^{(1)}_{21} = \mtx{0}$ and  $\| \mtx{A}^{(1)}_{12} \|$ is small. At this point, we may estimate the norm by finding the norms of the small $\mtx{A}^{(1)}_{11}$ block and the large $\mtx{A}^{(1)}_{22}$ block separately. $\| \mtx{A}^{(1)}_{11} \|_*$ may be computed efficiently with an SVD computation, and we may recursively apply our strategy for $\mtx{A}$ to $\mtx{A}^{(1)}_{22}$ to estimate its norm.

% -----------------------------------------------------------------------------
\subsection{Computing $\mtx{V}$}
\label{sec:compute_V}
% -----------------------------------------------------------------------------

At this point in the derivation, the most pressing question is how to find $\mtx{U}$ and $\mtx{V}$ that satisfy the aforementioned conditions. We first consider the construction of $\mtx{V}$, at which point $\mtx{U}$ will easily follow. Recent work in randomized subspace iteration, discussed in papers including \cite{halko2011finding,rokhlin2009randomized, woolfe2008fast,halko2011algorithm}, enables the efficient computation of the desired $\mtx{V}$.

In particular, consider a sampling matrix $\mtx{Y} = \mtx{A}^* \mtx{G}$, where $\mtx{G}$ is an $m \times b$ Gaussian random matrix. Then with high probability, the column space of  $\mtx{Y}$ often aligns closely to the subspace spanned by the leading $b$ right singular vectors of $\mtx{A}$. As discussed in \cite{rokhlin2009randomized,halko2011finding}, in certain situations when the decay of singular values of $\mtx{A}$ makes the alignment suboptimal, then using $\mtx{Y} = (\mtx{A}^* \mtx{A})^q \mtx{A}^* \mtx{G}$ instead provides considerable correction, where $q$ is some small nonnegative integer. Typically, $q = 0,1,2$ suffices to obtain a close-to-optimal approximation to the desired subspace.

Thus, the following steps compute a matrix $\mtx{V}$ whose columns span approximately the same subspace as the leading $b$ right singular vectors of $\mtx{A}$:

\begin{itemize}
\item Draw an $n \times b $ random Gaussian matrix $\mtx{G}$.
\item Compute $\mtx{Y}  = (\mtx{A}^* \mtx{A})^q \mtx{A}^* \mtx{G} $ for some small nonnegative integer $q$.
\item Compute an orthonormal basis for the column space of $\mtx{Y}$ with a QR factorization using Householder reflectors to obtain $\mtx{Y} = \mtx{V} \mtx{R}$.
\end{itemize}

% -----------------------------------------------------------------------------
\subsection{Computing $\mtx{U}$}
% -----------------------------------------------------------------------------

To compute $\mtx{U}$, we first note that if our computed $\mtx{V}$ were optimal, then the first $b$ columns of $\mtx{A} \mtx{V}$ would be composed exclusively of linear combinations of the leading $b$ left singular vectors of $\mtx{A}$. Thus, to form $\mtx{U}$, we may simply compute an orthonormal basis for the first $b$ columns of $\mtx{A} \mtx{V}$ via Householder reflectors to obtain $\mtx{A} \mtx{V}(:,1:b) = \mtx{U} \mtx{R}$.

\begin{remark}
In applications, the input matrix $\mtx{A}$ is often rank-deficient, with numerical rank $k$ determined by the $k$ for which $\sigma_i(\mtx{A})$ is less than some specified tolerance for $k < i \le n$. In the case where $k \ll n$, \texttt{randNN} may be sped up substantially by terminating the algorithm once the estimated singular values drop below some user-defined threshold.
\end{remark}

\begin{figure}

\begin{tabular}{|c|c|}\hline
\begin{minipage}{0.52\textwidth}
\small
\begin{verbatim}

function ss = randNN(A,b,q)
  m = size(A,1); n = size(A,2);
  ss = zeros(n,1);
  T = A;
  for i=1:ceil(n/b)
    I1 = (b*(i-1)+1):min((b*i),n);
    I2 = (b*i+1):m;
    J1 = (b*(i-1)+1):min((b*i),n);
    J2 = (b*i+1):n;
    if isempty(J2) == 0
      T_work = T([I1 I2],[J1 J2]);
      [TT,ss_part] = step_nn(T_work,b,q);
      ss(I1) = ss_part;
      T([I1 I2],[J1 J2]) = TT;
    else
      ss_part = svd(T(I1,J1));
      ss(I1) = ss_part;
    end
  end
return

\end{verbatim}
\end{minipage}
&
\begin{minipage}{0.45\textwidth}
\small
\begin{verbatim}

function [T,ss_part] = step_nn(A,b,q)
  G = randn(size(A,1),b);
  Y = A'*G;
  for i=1:q
    Y = A'*(A*Y);
  end
  [V,~] = qr(Y);
  T = A*V;
  [U,R] = qr(T(:,1:b));
  T(:,1:b) = R;
  T(:,(b+1):end) = U'*T(:,(b+1):end);
  ss_part = svd(R);
return








\end{verbatim}
\end{minipage}
\\[2mm] \hline
\end{tabular}
\caption{Matlab code for the algorithm \texttt{randNN} that, given an $m\times n$ matrix $\mtx{A}$,
computes estimates for each of its singular values.
The input parameters $b$ and $q$ reflect the block size and the number of steps of power iteration,
respectively. This code is simplistic in that it does not store or apply the transformation matrices $\mtx{U}$ and $\mtx{V}$ efficiently as products of Householder
vectors. It also does not apply a stopping criterion to terminate the algorithm if the estimated singular values become small, nor does it compute the error bound discussed in Section \ref{sec:err_bds}.}
\label{fig:randNNmatlab}
\end{figure}

% =============================================================================
\section{Error bounds}
\label{sec:err_bds}
% =============================================================================

In this section, we make note of an error bound for the accuracy of the estimated singular values resulting from the \texttt{randNN} algorithm.

First, note that after applying $\lceil n/b \rceil$ left and right rotation matrices to $\mtx{A}$ as described in \ref{sec:alg}, we are left with an upper triangular matrix which we shall call $\mtx{T}$ whose mass is concentrated in blocks of size $b \times b$ along its diagonal. Therefore, let $\mtx{T}_d$ be the block diagonal matrix consisting of the $b \times b$ main diagonal blocks of $\mtx{T}$, and let $\mtx{T}_u$ be the upper triangular matrix defined by the relation
\[
\mtx{T} = \mtx{T}_d + \mtx{T}_u.
\]
Observe that the estimated singular values of $\mtx{A}$ in our algorithm consist precisely of the computed singular values of $\mtx{T}_d$. Next, since $\mtx{A}$ and $\mtx{T}$ only differ by orthogonal transforms, we have that $\| \mtx{A} \|_* = \| \mtx{T} \|_*$, and more specifically, $\sigma_i(\mtx{A}) = \sigma_i(\mtx{T})$ for $i=1,2,\ldots,\min(m,n)$. Thus an error bound given in \cite{stewart1998perturbation} and \cite{mirsky1960symmetric} states that
\begin{equation}
\sqrt{\sum_{i=1}^{\min(m,n)} \bigl(\sigma_i(\mtx{A}) - \sigma_i(\mtx{T}_d)\bigr)^2} \le \| \mtx{T}_u \|_F.
\end{equation}

Since the calculation of $\| \mtx{T}_u \|_F$ is $\mathcal{O}(n^2)$, this bound may be calculated at negligible relative cost as part of \texttt{randNN}, serving as an assurance of the estimation's validity.

% =============================================================================
\section{Numerical experiments}
\label{sec:numerics}
% =============================================================================

% -----------------------------------------------------------------------------
\subsection{Computational speed}
% -----------------------------------------------------------------------------

In this section, we investigate the speed of the proposed algorithm \texttt{randNN} and compare it to the speed of an exact computation of the singular values using LAPACK's highly optimized \texttt{dgesvd}.

All experiments reported in this note were performed on an Intel Core i7-6700K processor (4.0 GHz) with 4 cores. In order to be able to show scalability results, the clock speed was throttled at 4.0 GHz, turning off so-called turbo boost. Other details of interest include that the OS used was Ubuntu (Version 16.04.2), and the code was compiled with Intel's \texttt{icc} (Version 17.0.4.196). Main routines for computing the singular values (\texttt{dgesvd}) were taken from Intel's MKL library (Version 2017.0.3). Standard BLAS routines used in our implementation of \texttt{randNN} were also taken from the Intel MKL library.

Each of the three algorithms we tested (\texttt{randNN},\texttt{randUTV},\texttt{SVD}) was applied to double-precision real matrices of size $n \times n$. We report the following times:
\begin{itemize}
\item $T_{\text{svd}}$. The time in seconds for the LAPACK function \texttt{dgesvd} from Intel's MKL, where the orthogonal matrices $\mtx{U}$ and $\mtx{V}$ in the SVD were not built.
\item $T_{\text{randUTV}}$. The time in seconds for the function \texttt{randUTV} described in \cite{martinsson2017randutv}. The authors' original implementation was used, and the orthogonal matrices $\mtx{U}$ and $\mtx{V}$ in the UTV decomposition were not built.
\item $T_{\text{randNN}}$. The time in seconds for our implementation of \texttt{randNN}.
\end{itemize}

In all cases, we used a block size of $b=64$. While likely not optimal for all problem sizes, this block size yields near best performance and allows us to easily compare and contrast the performance of the different implementations.

\begin{figure}
\includegraphics[width=.45\textwidth]{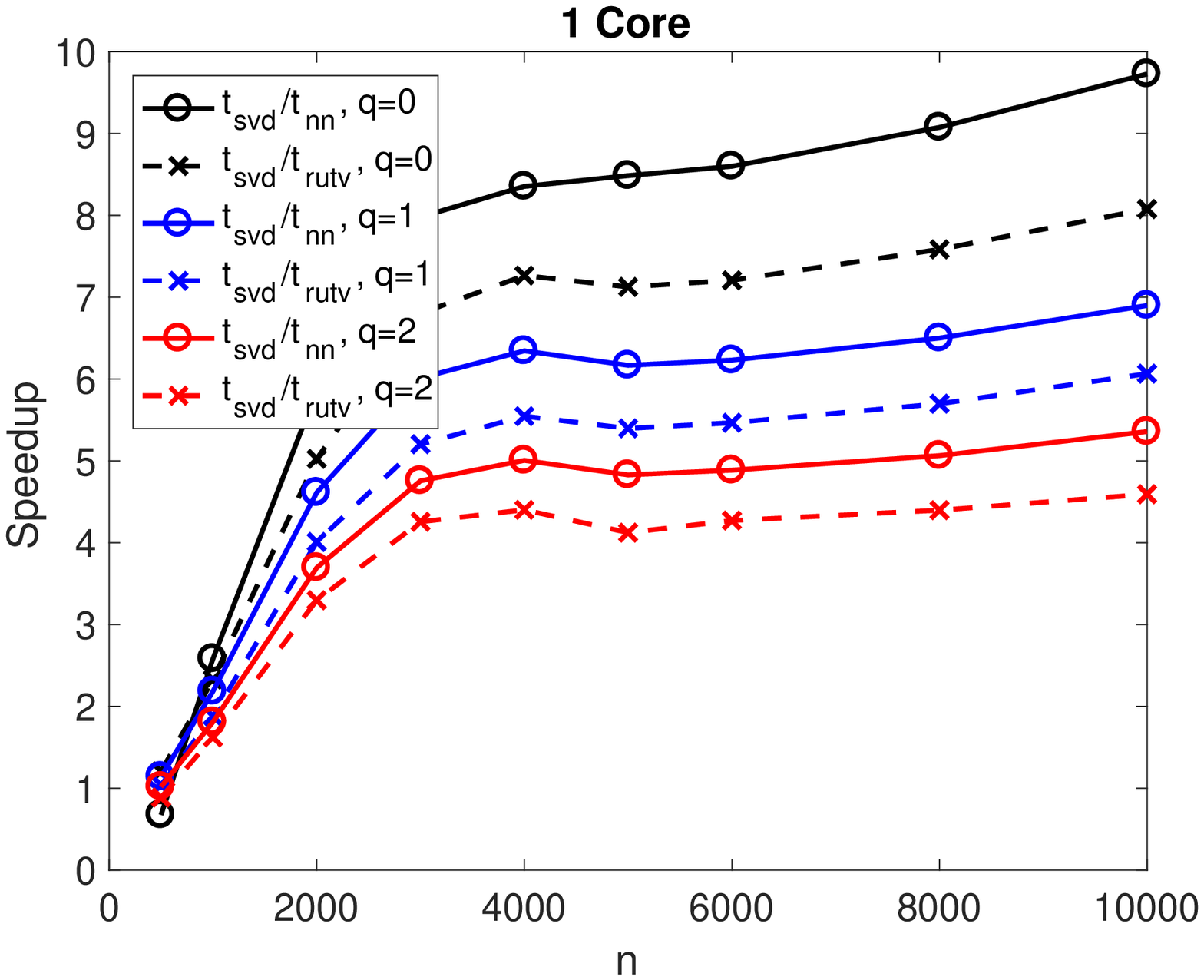}
\includegraphics[width=.45\textwidth]{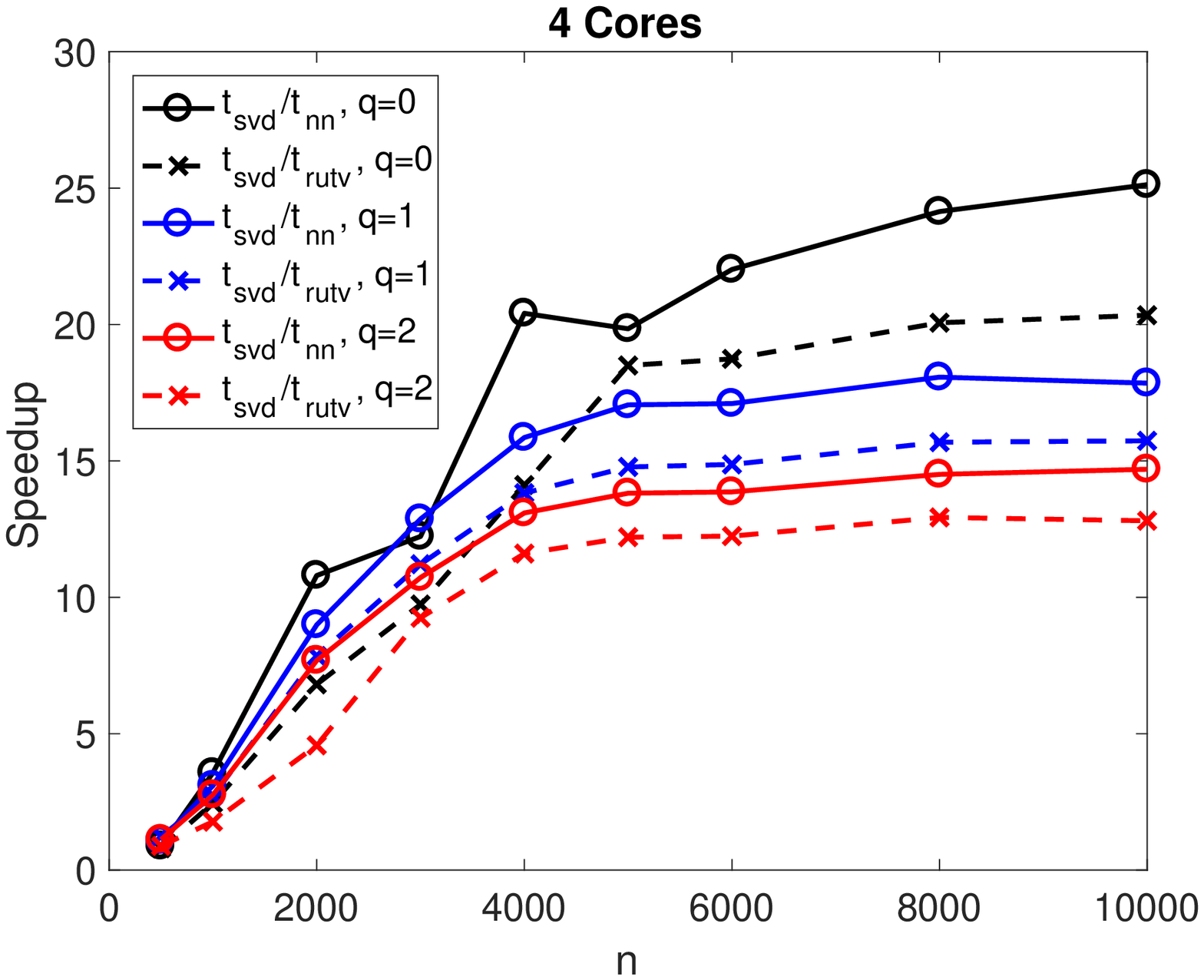}
\caption{Computational speed of \texttt{randNN} compared to the speeds of the LAPACK routine \texttt{dgesvd} and the recently introduced \texttt{randUTV} \cite{martinsson2017randutv}. The algorithms were applied to double-precision real matrices of size $n \times n$. \label{fig:speed}}
\end{figure}

As expected, Figure \ref{fig:speed} shows that both \texttt{randUTV} and \texttt{randNN} are markedly faster than LAPACK's \texttt{dgesvd}, even in the single case. When more cores are added to the computation, the gain in speed is even more dramatic, hi-lighting the fact that the SVD algorithm is not designed to make efficient use of parallel computing architectures. Finally the increase in speed from \texttt{randUTV} to \texttt{randNN} is modest, but in a line search setting where the nuclear norm must be computed many times, this saving can add up over the course of the solution of the rank minimization problem.

% -----------------------------------------------------------------------------
\subsection{Errors}
\label{sec:error}
% -----------------------------------------------------------------------------

\begin{figure}
\includegraphics[width=.45\textwidth]{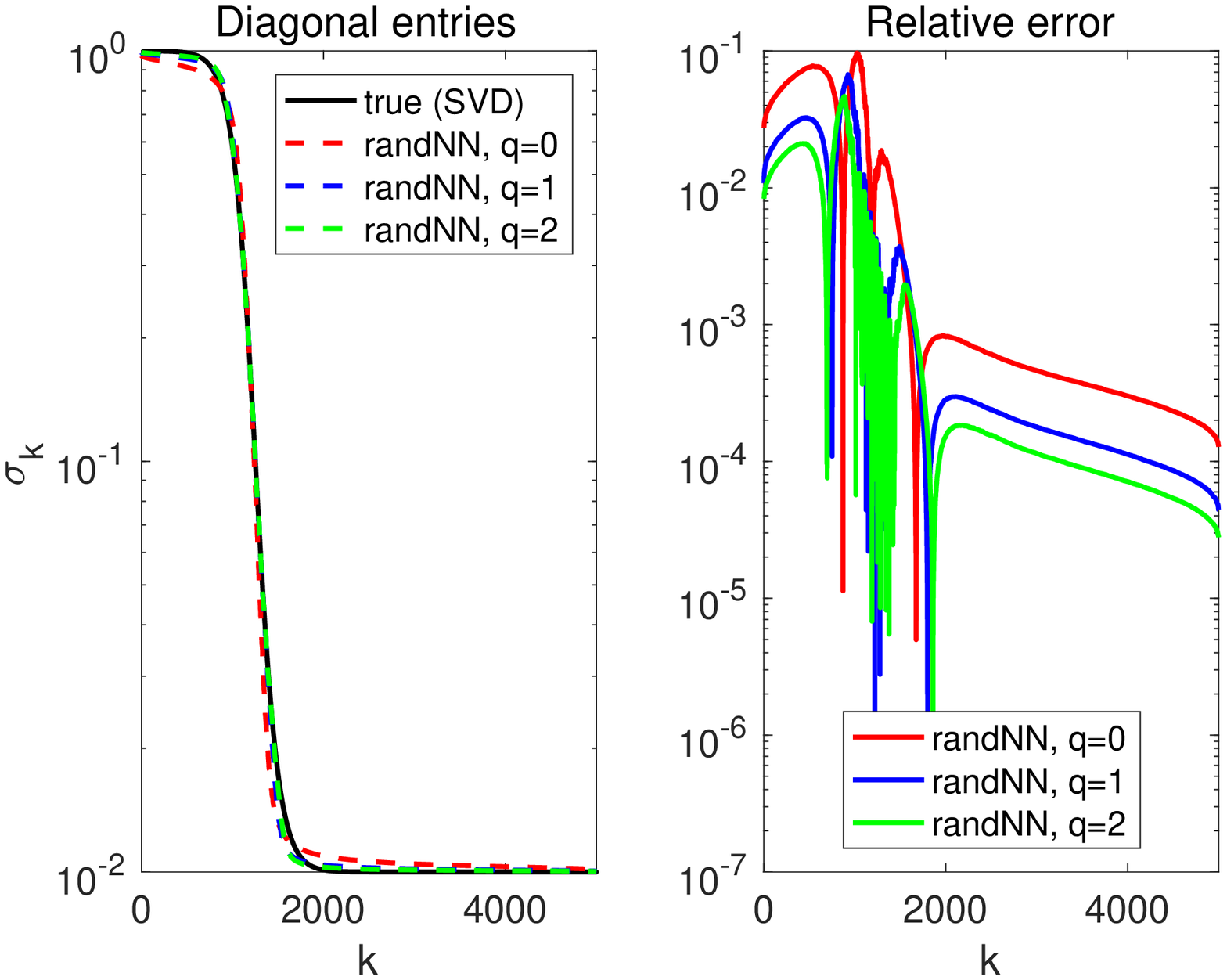}
\includegraphics[width=.45\textwidth]{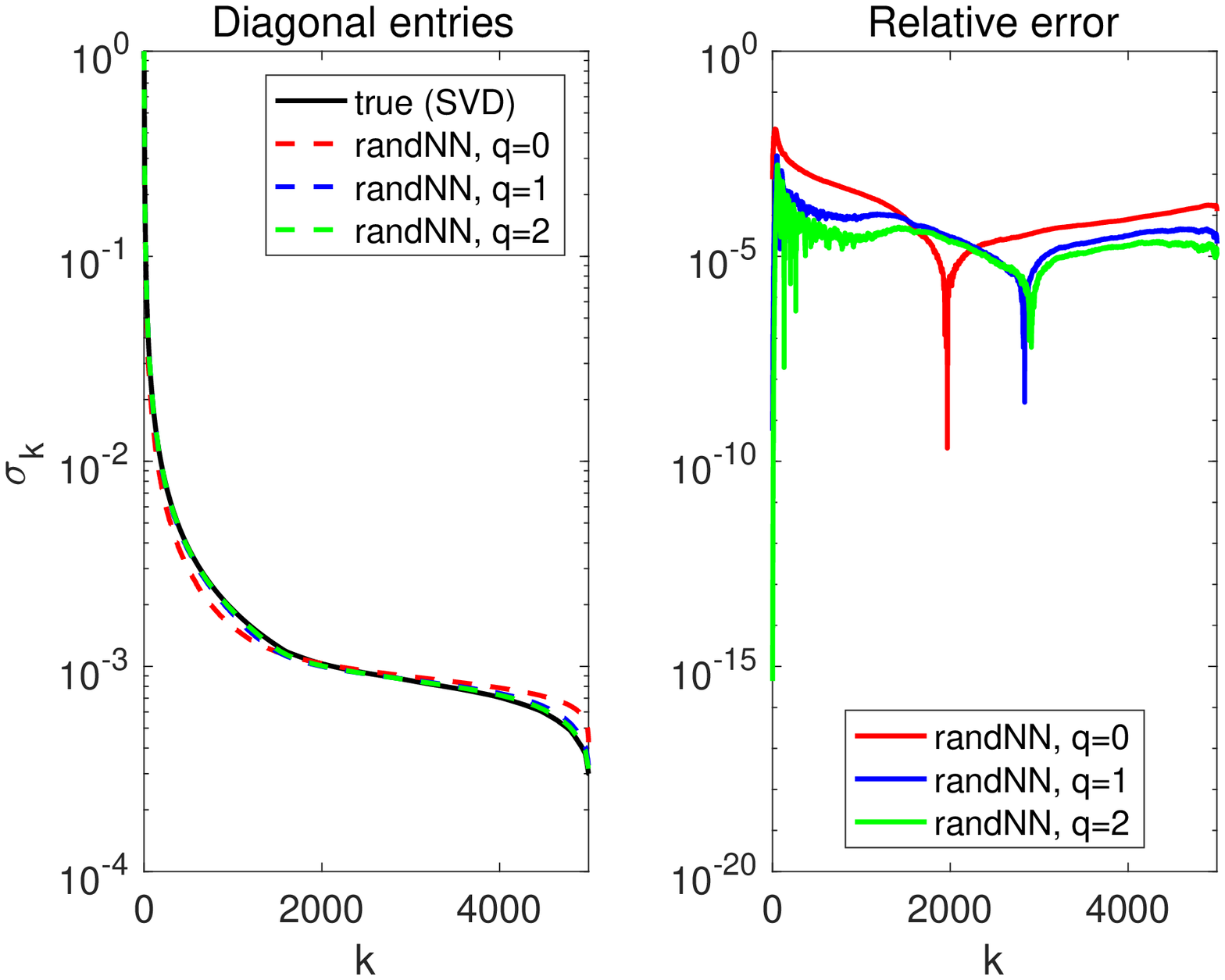} \\

\caption{Errors for the test matrices described in Section \ref{sec:error}. Left: ``Matrix 1,'' the matrix artificially to have the singular value decay pattern shown here. Right: ``Matrix 2,'' a matrix resulting from the discretization of of Boundary Integral Equation. \label{fig:error}}
\end{figure}

In this section, we report the results of numerical experiments that were conducted to test the accuracy of the approximation to the singular values provided by \texttt{randNN}. Specifically, we compare the estimated singular values to the true singular values of two different test matrices:

\begin{itemize}
\item \textit{Matrix 1 (S shaped decay):} This is an $n \times n$ matrix of the form $\mtx{A} = \mtx{U} \mtx{D} \mtx{V}^*$ where $\mtx{U}$ and $\mtx{V}$ are randomly drawn matrices with orthonormal columns (obtained by performing QR on a random Gaussian matrix), and where the diagonal entries of $\mtx{D}$ are chosen to first hover around 1, then decay rapidly, and then level out at $10^{-6}$, as shown in Figure \ref{fig:error} (black line) on the left.
\item \textit{Matrix 2 (Single Layer BIE):} This matrix is the result of discretizing a Boundary Integral Equation (BIE) defined on a smooth closed curve in the plane. To be precise, we discretized the so called ``single layer'' operator associated with the Laplace equation using a $6^{\rm th}$ order quadrature rule designed by Alpert \cite{alpert1999hybrid}. This is a well-known ill-conditioned problem for which column pivoting is essential in order to stably solve the corresponding linear system.
\end{itemize}

In each case, a matrix size of $m = n = 5000$ was used with a ``power iteration parameter'' (see Section \ref{sec:compute_V}) of $q = 2$ and a block size of $b = 64$.

For each test matrix, we plot both the estimated and true singular values themselves and the relative error of the estimates
\[
\frac{| \sigma_i(\mtx{A}) - \sigma_i(\mtx{T}_d) |}{| \sigma_i(\mtx{A}) |}, \quad i=1,2,\ldots,\min(m,n).
\]

% =============================================================================
\section{Availability of code}
% =============================================================================

An implementation of the discussed algorithm is available under 3-clause (modified) BSD license from:
\[
\text{https://github.com/nheavner/nn\_code}
\]

\bibliography{main_bib}
\bibliographystyle{amsplain}

\end{document}